\documentclass[11pt]{article}
\usepackage{amsmath}
\usepackage{graphicx}
\def\disp{\displaystyle}

\def\Limsup{\mathop{{\rm Lim}\,{\rm sup}}}

\def\tto{\;{\lower 1pt \hbox{$\rightarrow$}}\kern -10pt
\hbox{\raise 2pt \hbox{$\rightarrow$}}\;}
\def\Hat{\widehat}

\def\Bar{\overline}
\def\ra{\rangle}
\def\la{\langle}
\def\ve{\varepsilon}
\def\B{I\!\!B}
\def\h{\hfill\Box}
\def\R{I\!\!R}
\def\N{I\!\!N}
\def\ox{\bar{x}}

\def\dom{\mbox{\rm dom}\,}

\def\cone{\mbox{\rm cone}\,}

\def\i{\mbox{\rm int}\,}
\def\substack#1#2{{\scriptstyle{#1}\atop\scriptstyle{#2}}}

\def\h{\hfill\triangle}
\def\dn{\downarrow}
\def\O{\Omega}
\def\ph{\varphi}
\def\emp{\emptyset}
\def\st{\stackrel}
\def\oR{\Bar{\R}}

\def\dd{\delta}
\def\al{\alpha}

\def \N{I\!\!N}

\newcounter{lk}

\begin{document}
\begin{center}
\vspace*{0.3in} {\bf APPLICATIONS OF VARIATIONAL ANALYSIS TO A GENERALIZED FERMAT-TORRICELLI PROBLEM}\\[3ex]
BORIS S. MORDUKHOVICH\footnote{Department of Mathematics, Wayne
State University, Detroit, MI 48202, USA (email:
boris@math.wayne.edu). Research of this author was partially
supported by the US National Science Foundation under grants
DMS-0603846 and DMS-1007132 and by the Australian Research Council
under grant DP-12092508} and NGUYEN MAU NAM\footnote{Department of
Mathematics, The University of Texas--Pan American, Edinburg,
TX 78539--2999, USA (email: nguyenmn@utpa.edu).}\\[2ex]
\end{center}
\small{\bf Abstract.} In this paper we develop new applications of
variational analysis and generalized differentiation to the
following optimization problem and its specifications: given $n$
closed subsets of a Banach space, find such a point for which the
sum of its distances to these sets is minimal. This problem can be
viewed as an extension of the celebrated Fermat-Torricelli problem:
given three points on the plane, find another point such that the
sum of its distances to the designated points is minimal. The
generalized Fermat-Torricelli problem formulated and studied in this
paper is of undoubted mathematical interest and is promising for
various applications including those frequently arising in location
science, optimal networks, etc. Based on advanced tools and recent
results of variational analysis and generalized differentiation, we
derive necessary as well as necessary and sufficient optimality
conditions for the extended version of the Fermat-Torricelli problem
under consideration, which allow us to completely solve it in some
important settings. Furthermore, we develop and justify a numerical
algorithm of the subgradient type to find optimal solutions in convex settings and provide its numerical implementations.\\[2ex]
\newtheorem{Theorem}{Theorem}[section]
\newtheorem{Proposition}[Theorem]{Proposition}
\newtheorem{Remark}[Theorem]{Remark}
\newtheorem{Lemma}[Theorem]{Lemma}

\newtheorem{Corollary}[Theorem]{Corollary}
\newtheorem{Definition}[Theorem]{Definition}
\newtheorem{Example}[Theorem]{Example}
\renewcommand{\theequation}{\thesection.\arabic{equation}}
\normalsize

\section{Introduction and Problem Formulation}
\setcounter{equation}{0}

In the early 17th century Pierre de Fermat proposed the following
problem: given three points on the plane, find a fourth point such
that the sum of its Euclidean distances to the three given points is
minimal. This problem was solved by Evangelista Torricelli and was
named the {\em Fermat-Torricelli problem}. Torricelli's solution
states the following: if none of the interior angles of the triangle
formed by the three fixed points reaches or exceeds $120^\circ$, the
minimizing point in question is located inside this triangle in such
a way that each side of the triangle is seen at an angle of
$120^\circ$; otherwise it is the obtuse vertex of the triangle. This
point is often called the Fermat-Torricelli point.

In the 19th century Jakob Steiner examined this problem in further
depth and extended it to include a finitely many points on the
plane. A number of other extensions have been proposed and studied
over the years. This and related topics have nowadays attracted
strong attention of many mathematicians and applied scientists; see,
e.g., \cite{b,bms,msw,tan} with the references therein for the
history, various extensions, modifications, and applications to
location science, statistics, optimal networks, etc. Note that,
despite beautiful solutions obtained for particular extensions of
the Fermat-Torricelli problem, we are not familiar with theoretical
methods and/or numerical algorithms developed in rather general
settings. We particularly refer the reader to \cite{bz,k,wp} and the
bibliographies therein to Weiszfeld's algorithm and its
modifications for the problem of minimizing weighted sums of
Euclidean norms (also known as the Weber problem) and to \cite{acco}
for efficient interior point type methods for similar problems in
finite-dimensional spaces. \vspace*{0.05in}

In this paper we study a far-going generalization of the
Fermat-Torricelli problem that is formulated below. It extends, in
particular, a generalized version of the classical Steiner (and
Weber) versions with replacing given points therein by a finitely
many {\em closed sets} in Banach spaces. Furthermore, our new
extension of the Fermat-Torricelli problem covers a fast majority of
the previous ones and seems to be interesting for both the theory
and applications to various location models, optimal networks,
wireless communications, etc.

We propose to employ powerful tools of modern variational analysis
and generalized differentiation to study the extended version of
Fermat-Torricelli problem and its specification from the
theoretical/qualitative and numerical/algorithmic viewpoints. In the
first direction our goal is to derive necessary as well as necessary
and sufficient optimality conditions for generalized
Fermat-Torricelli points and then to use them for explicit
determining these points in some remarkable settings. Our numerical
analysis involves developing an algorithm of the subgradient type
and considering its specifications and implementations in the case
of the generalized Fermat-Torricelli problem determined by an
arbitrary number of convex sets in finite-dimensional spaces.
\vspace*{0.05in}

Let us now formulate the generalized Fermat-Torricelli problem of
our study. Consider the so-called {\em minimal time function}
\begin{equation}\label{minimal time}
T^F_\O(x):=\inf\big\{t\ge 0\big|\;\O\cap (x+tF)\ne\emp\big\}
\end{equation}
with the {\em constant dynamics} $\dot x\in F$ described by a
closed, bounded, and convex subset $F\ne\emp$ of a Banach space $X$
and with the closed {\em target} set $\O\ne\emp$ in $X$; these are
our {\em standing assumptions} in this paper. We refer the reader to
\cite{cgm,cowo,heng,mn10} and the bibliographies therein for various
results on minimal time functions and their applications. When $F$
is the closed unit ball $\B$ of $X$, the minimal time function
(\ref{minimal time}) becomes the standard {\em distance function}
\begin{equation}\label{distance function1}
\mbox{d}(x;\O)=\inf\big\{\|x-\omega\|\;\big|\;\omega\in\O\big\}
\end{equation}
generated by the norm $\|\cdot\|$ on $X$. Given now an arbitrary
number of closed subsets $\O_i\ne\emp$, $i=1,\ldots,n$, of $X$, we
introduce the {\em generalized Fermat-Torricelli problem} as
follows:
\begin{equation}\label{ft}
\mbox{minimize }\;T(x):=\sum_{i=1}^n T^F_{\O_i}(x),\quad x\in X.
\end{equation}
For $F=\B$ in \eqref{ft} this problem reduces to
\begin{equation}\label{distance function}
\mbox{minimize }\;\mbox{D}(x):=\sum_{i=1}^n\mbox{d}(x;\O_i),
\end{equation}
which corresponds to the Steiner-type extension of the
Fermat-Torricelli problem in Banach spaces when all the sets $\O_i$,
$i=1,\ldots,n$, are singletons. Observe that even in the latter
classical case the optimization problem \eqref{ft} and its
specification \eqref{distance function} are {\em nonsmooth} while
being convex if all the sets $\O_i$ have this property. It is thus
natural to study these problems by means of advanced tools of
variational analysis and generalized differentiation.

The rest of the paper is organized as follows. In Section~2 we
define and discuss basic tools of variational analysis needed for
formulations and proofs of the main results of this paper. Section~3
is devoted to computing and estimating subdifferentials of minimal
time functions that play a crucial role in our study of the
generalized Fermat-Torricelli problem and its specifications. In
Section~4 we derive necessary conditions for the general problem
\eqref{ft} in Banach spaces as well as necessary and sufficient
conditions in the case of convexity, which are then used for
complete descriptions of Fermat-Torricelli points in some important
settings. Finally, Section~5 presents and justifies a numerical
algorithm for solving the generalized Fermat-Torricelli problem with
convex data in finite-dimensional spaces.\vspace*{0.05in}

Throughout the paper we use standard notation and terminology of
variational analysis; see, e.g., \cite{mor06a,rw}. Recall that,
given a set-valued mapping $G\colon X\tto X^*$ between a Banach
space $X$ and its topological dual $X^*$, the {\em sequential
Painlev\'e-Kuratowski upper/outer limit} as $x\to\ox$ is defined by
\begin{eqnarray}\label{pk}
\begin{array}{ll}
\disp\Limsup_{x\to\ox}G(x):=\Big\{x^*\in X^*\Big|&\exists\,\mbox{
sequences } \;x_k\to\ox,\;x^*_k\st{w^*}{\to}x^*\;\mbox{ as }
\;k\to\infty\\
&\mbox{such that }\;x^*_k\in G(x_k)\;\mbox{ for all
}\;k\in\N:=\{1,2,\ldots\}\Big\},
\end{array}
\end{eqnarray}
where $w^*$ signifies the weak$^*$ topology of $X^*$. For a set
$\O\subset X$ the symbol $x\st{\O}{\to}\ox$ means that $x\to\ox$
with $x\in\O$. If $\ph\colon X\to\oR:=(-\infty,\infty]$ is an
extended-real-valued function finite at $\ox$, the symbol
$x\st{\ph}{\to}\ox$ signifies the convergence $x\to\ox$ with
$\ph(x)\to\ph(\ox)$.

\section{Tools of Variational Analysis}
\setcounter{equation}{0}

In this section we briefly review some basic constructions and
results of the generalized differentiation theory in variational
analysis  that are widely used in what follows. The reader can find
all the proofs, discussions, and additional material in the books
\cite{bz,mor06a,mor06b,rw,s} and the references therein in both
finite and infinite dimensions. Unless otherwise stated, all the
spaces under consideration are Banach with the norm $\|\cdot\|$ and
the canonical pairing $\la\cdot,\cdot\ra$ between the space in
question and its topological dual.

Let us start with {\em convex} functions $\ph\colon X\to\oR$. Given
$\ox\in\dom\ph:=\{x\in X|\;\ph(x)<\infty\}$, the {\em
subdifferential} (collection of subgradients) of $\ph$ at $\ox$ in
the sense of convex analysis is
\begin{equation}\label{convex subdifferential}
\partial\ph(\ox):=\big\{x^*\in X^*\big|\;\la x^*,x-\ox\ra\le\ph(x)-\ph(\ox)
\;\mbox{ for all }\;x\in X\big\}.
\end{equation}
Directly from definition \eqref{convex subdifferential} we have the
following nonsmooth counterpart of the classical Fermat stationary
rule for convex functions:
\begin{equation}\label{fermat}
\ox\;\mbox{ is a minimizer of }\;\ph\;\mbox{ if and only if
}\;0\in\partial\ph(\ox).
\end{equation}
The subdifferential of convex analysis \eqref{convex
subdifferential} satisfies a number of important calculus rules that
are mainly based on separation theorems for convex sets. The central
calculus result is the following Moreau-Rockafellar theorem for
representing the subdifferential of sums.

\begin{Theorem}\label{sum rule} {\bf (subdifferential sum rule for convex functions).}
Let $\ph_i\colon X\to\oR$, $i=1,\ldots,m$, be convex lower
semicontinuous functions on a Banach space $X$. Assume that there is
a point $\ox\in\cap_{i=1}^n\dom\ph_i$ at which all (except possibly
one) of the functions $\ph_1,\ldots,\ph_m$ are continuous. Then we
have the equality
\begin{equation*}
\partial\Big(\sum_{i=1}^m\ph_i\Big)(\ox)=\sum_{i=1}^m\partial\ph_i(\ox).
\end{equation*}
\end{Theorem}

Given a convex set $\O\subset X$ and a point $\ox\in\O$, the
corresponding geometric counterpart of \eqref{convex
subdifferential} is the {\em normal cone} to $\O$ at $\ox$ defined
by
\begin{equation}\label{cnc}
N(\ox;\O)=\big\{x^*\in X^*\big|\;\la x^*,x-\ox\ra\le 0\;\mbox{ for
all }\;x\in\O\big\},
\end{equation}
which is in fact the subdifferential \eqref{convex subdifferential}
of the set indicator function  $\dd(x;\O)$ at $\ox$ that is equal to
$0$ for $x\in\O$ and to $\infty$ for $x\notin\O$. \vspace*{0.05in}

Besides the aforementioned convex constructions suitable for the
study of the generalized Fermat-Torricelli problem \eqref{ft} in the
case of convex sets $\O_i$, we need in what follows their extensions
to nonconvex objects.

Given an arbitrary extended-real-valued function $\ph\colon X\to\oR$
finite at $\ox$ and given $\ve\ge 0$, define first the $\ve$-{\em
subdifferential} of $\ph$ at $\ox$ by
\begin{eqnarray}\label{2.3}
\Hat\partial_\ve\ph(\ox):=\Big\{x^*\in
X^*\Big|\;\liminf_{x\to\ox}\frac{\ph(x)-\ph(\ox)- \la
x^*,x-\ox\ra}{\|x-\ox\|}\ge-\ve\Big\}.
\end{eqnarray}
For $\ve=0$ the set $\Hat\partial\ph(\ox):=\Hat\partial_0\ph(\ox)$
is known as {\em regular/viscosity/Fr\'echet subdifferential} of
$\ph$ at $\ox$; it reduces to the classical gradient
$\{\nabla\ph(\ox)\}$ when $\ph$ is Fr\'echet differentiable at this
point and to the subdifferential \eqref{convex subdifferential} when
$\ph$ is convex. However, being naturally and rather simply defined,
the Fr\'echet subdifferential and its $\ve$-enlargements \eqref{2.3}
do not possess---apart from locally convex settings and the like---a
number of required calculus and related properties. For example,
$\Hat\partial\ph(\ox)$ may often be empty (e.g., for $\ph(x)=-|x|$)
and an analog of the sum rule from Theorem~\ref{convex
subdifferential} does not hold for $\Hat\partial\ph(\ox)$ whenever
$\ve\ge 0$; e.g., in the case of $\ph_1(x)=|x|$ and $\ph_2(x)=-|x|$.

The situation dramatically changes when we employ a sequential
regularization of the $\ve$-subdifferentials \eqref{2.3} defined by
\begin{eqnarray}\label{2.4}
\partial\ph(\bar x):=\Limsup_\substack{x\xrightarrow{\ph}\bar x}{\ve\dn 0}
\Hat\partial_\ve\ph(x)
\end{eqnarray}
via the sequential outer limit \eqref{pk} and known as the {\em
basic/limiting/Mordukhovich subdifferential} of $\ph$ at $\ox$. We
can equivalently put $\ve=0$ in \eqref{2.4} if $\ph$ is lower
semicontinuous around $\ox$ and if the space $X$ is {\em Asplund},
i.e., each of its separable subspaces has a separable dual. The
latter subclass of Banach spaces is sufficiently large including, in
particular, every reflexive space and every space with a separable
dual. On the other hand, it does not contain some classical Banach
spaces important for applications as, e.g., $C[0,1]$ and $L^1[0,1]$.

A geometric counterpart of the subdifferential \eqref{2.4} is the
corresponding (basic, limiting, Mordukhovich) {\em normal cone} to a
set $\O\subset X$ at $\ox\in\O$ that can be defined via the
subdifferential \eqref{2.4} of the indicator function
 $N(\ox;\O):=\partial\dd(\ox;\O)$ and reduces to the normal cone of
convex analysis \eqref{cnc} for convex sets $\O$. The given
definition of our basic normal can be equivalently rewritten in the
limiting form
\begin{eqnarray}\label{2.3.1}
N(\ox;\O)=\Limsup_\substack{x\to\ox}{\ve\dn 0}\Hat N_\ve(x;\O)
\end{eqnarray}
with the sets of $\ve$-normals $\Hat N_\ve(\cdot;\O)$ defined for
$\ve\ge 0$ by
\begin{eqnarray}\label{1.2.1}
\Hat N_\ve(\ox;\O):=\disp\Big\{x^*\in
X^*\Big|\;\limsup_{x\st{\O}{\to}\ox}\frac{\la
x^*,x-\ox\ra}{\|x-\ox\|}\le\ve\Big\},\quad\ox\in\O,
\end{eqnarray}
where $\Hat N_\ve(\ox;\O):=\emp$ if $\ox\notin\O$ for convenience.
When the set $\O$ is locally closed around $\ox$ and the space $X$
is Asplund, we can equivalently replace $\Hat N_\ve(\cdot;\O)$ in
\eqref{1.2.1} by the {\em prenormal/Fr\'echet normal cone} $\Hat
N(\cdot;\O):=\Hat N_0(\cdot;\O)$. Furthermore, in the case of
$X=\R^n$ the normal cone \eqref{2.3.1} admits the representation
\begin{eqnarray}\label{nc}
N(\ox;\O)=\disp\Limsup_{x\to\ox}\big[\mbox{cone}\big(x-\Pi(x;\O)\big)\big],
\end{eqnarray}
where $\Pi(x;\O)$ denotes the Euclidean projection of the point
$x\in\R^n$ onto the closed set $\O$, and where $\cone\O$ signifies
the collection of rays spanned on $\O$. Representation \eqref{nc}
was actually the original definition of the limiting normal cone in
\cite{mor76}.

In spite of the nonconvexity of the limiting constructions
\eqref{2.4} and \eqref{2.3.1}, they enjoy well-developed calculus
rules that are pretty comprehensive in the Asplund space setting and
are based on variational/extremal principles; see, e.g.,
\cite{mor06a}. In particular, the following sum rule for the
subdifferential \eqref{2.4} is used in this paper.

\begin{Theorem}\label{nonconvex sum rule} {\bf (subdifferential sum rule
for nonconvex functions).} Let $\ph_i\colon X\to\oR$,
$i=1,\ldots,m$, be lower semicontinuous functions on an Asplund
space $X$. Suppose that all (except possibly one) of them are
locally Lipschitzian around $\ox\in\cap_{i=1}^m\dom\ph_i$. Then we
have the inclusion
\begin{equation*}
\partial\Big(\sum_{i=1}^n\ph_i\Big)(\ox)\subset\sum_{i=1}^n\partial\ph_i(\ox).
\end{equation*}
\end{Theorem}

\section{Generalized Differentiation of Minimal Time Functions}
\setcounter{equation}{0}

This section is devoted to reviewing, for the reader's convenience,
some recent results on generalized differentiation of the minimal
time functions of type \eqref{minimal time} developed in detail in
our separate paper \cite{bmn10}, which in fact was mainly motivated
by the application to the generalized Fermat-Torricelli problem
\eqref{ft} given in what follows. Let us present and discuss the
major required results on generalized differentiation of the minimal
time functions in both convex and nonconvex cases. We say that
$\ox\in X$ is a {\em in-set point} for the minimal time function
\eqref{minimal time} if $\ox\in\O$ and that $\ox$ is an {\em
out-of-set point} for \eqref{minimal time} if $\ox\notin\O$.

The following result, which is a consequence of Proposition~4.1,
Proposition~4.2, and Theorem~5.2 from \cite{bmn10}, provides precise
relationships between basic subdifferential \eqref{2.4} of the
minimal time function \eqref{minimal time} and the basic normal cone
\eqref{nc} to the corresponding targets in the case of {\em in-set}
points.

\begin{Theorem}\label{lim-min} {\bf (basic subgradients of minimal
time functions and basic normals to targets at in-set points).} Let
$\ox\in\O$ for the minimal time function \eqref{minimal time} on a
Banach space $X$, and let the support level set $C^*\subset X^*$ be
defined by
\begin{equation}\label{sl}
C^*:=\big\{x^*\in X^*\big|\;\sigma_F(-x^*)\le 1\big\}
\end{equation}
via the support function of the dynamics given by
\begin{equation}\label{sf}
\sigma_F(x^*):=\sup_{x\in F}\la x^*,x\ra,\quad x^*\in X^*.
\end{equation}
Then we have the subdifferential upper estimate
\begin{eqnarray}\label{limiting1}
\partial T^F_\O(\ox)\subset N(\ox;\O)\cap C^*.
\end{eqnarray}
Furthermore, the latter holds as the equality
\begin{eqnarray}\label{limiting2}
\partial T^F_\O(\ox)=N(\ox;\O)\cap C^*
\end{eqnarray}
provided that the target set $\O$ is convex.
\end{Theorem}

The next result gives an upper estimate of the basic subdifferential
of the generally nonconvex minimal time function \eqref{minimal
time} at {\em out-of-set} points via basic subgradients of the
corresponding {\em Minkowski gauge}
\begin{eqnarray}\label{mk}
\rho_F(x):=\inf\big\{t\ge 0\big|\;x\in tF\big\},\quad x\in X,
\end{eqnarray}
associated with the dynamics $F$ and via basic normals to the target
set $\O$ at points belonging the the {\em minimal time/generalized
projection} of $\ox\notin\O$ to $\O$ defined by
\begin{eqnarray}\label{pr}
\Pi^F_\O(\ox):=\big(\ox+T^F_\O(\ox)F\big)\cap\O.
\end{eqnarray}
It is easy to see that the generalized projection \eqref{pr} reduces
to the standard (metric) one when $F=\B$, i.e., when \eqref{minimal
time} becomes the distance function \eqref{distance function1}.

To proceed, we recall that the minimal time function \eqref{minimal
time} is {\em well posed} at $\ox\notin\O$ with $T^F_\O(\ox)<\infty$
if for any sequence $x_k\to\ox$ with $T^F_\O(x_k)\rightarrow
T^F_\O(\ox)$ as $k\to\infty$ there is a sequence of projection
points $w_k\in\Pi^F_\O(x_k)$ containing a convergent subsequence.
This property is defined and discussed in \cite{bmn10}: cf.\ also
\cite[Subsection~1.3.3]{mor06a} for the case of distance functions.
The following conditions are sufficient for
well-posedness:\vspace*{0.05in}

$\bullet$ The target $\O$ is a compact subset of $X$;

$\bullet$ The space $X$ is finite-dimensional and $\O$ is a closed
subset of $X$;

$\bullet$ $X$ is reflexive, $\O\subset X$ is closed and convex, and
the Minkowski gauge associated with $F$ generates an
 equivalent {\em Kadec norm} on $X$, i.e., such that the weak and
norm convergences agree on the boundary of the unit sphere of
$X$.\vspace*{0.05in}

Here is the aforementioned result; cf.\ \cite[Theorem~6.3]{bmn10}.

\begin{Theorem}\label{mink3} {\bf (basic subgradients of minimal time
functions at out-of-set points via projections).} Let $\ox\notin\O$
with $T^F_\O(\ox)<\infty$, and let the minimal time function
\eqref{minimal time} be well posed at $\ox$. Then we have the upper
estimate
\begin{eqnarray}\label{wel1}
\partial T^F_\O(\ox)\disp\subset\bigcup_{\bar w\in\Pi^F_\O(\ox)}\big[-\partial
\rho_F(\bar w-\ox)\cap N(\bar w;\O)\big].
\end{eqnarray}
\end{Theorem}

Finally in this section, consider the case of {\em convexity} of the
minimal time function \eqref{minimal time}, which is equivalent to
the convexity of its target set $\O$ as shown, e.g., in
\cite[Proposition~3.6]{bmn10}. In this case we have some specific
results, which are not satisfied for general nonconvex minimal time
functions; see \cite{bmn10}  for more details. In particular, the
convex case allows us to establish important connections between the
basic subdifferential of \eqref{minimal time} and the corresponding
normal cone to the target {\em enlargements}
\begin{eqnarray}\label{en}
\O_r:=\big\{x\in X\big|\;T^F_\O(x)\le r\big\},\quad r>0,
\end{eqnarray}
at out-of-set points $\ox\notin\O$. The following result taken from
\cite[Theorem~7.3]{bmn10} contains what we need for applications to
the generalized Fermat-Torricelli problem in this paper.

\begin{Theorem}\label{convex case 2} {\bf (subgradients of convex minimal
time functions at out-of-set points).} Let the minimal time
 function \eqref{minimal time} be convex, and let $\ox\notin\O$ be
such that $\Pi^F_\O(\ox)\ne\emp$ with $r=T^F_\O(\ox)<\infty$ in
\eqref{en}. Then for any $\bar w\in\Pi^F_\O(\ox)$ we have the
relationships
\begin{eqnarray}\label{c7}
\begin{array}{ll}\partial T^F_\O(\ox)&=N(\ox;\O_r)\cap\big[-\partial
\rho_F(\bar w-\ox)\big]\\
&\subset N(\bar w;\O)\cap\big[-\partial\rho_F(\bar w-\ox)\big].
\end{array}
\end{eqnarray}
If in addition $0\in F$, then the inclusion in \eqref{c7} holds as
equality, and thus
\begin{eqnarray}\label{c8}
\partial T^F_\O(\ox)=N(\bar w;\O)\cap\big[-\partial\rho_F(\bar w-\ox)\big].
\end{eqnarray}
\end{Theorem}

\section{Generalized Fermat-Torricelli Problem: Optimality Conditions in
Finite and Infinite Dimensions}
\setcounter{equation}{0}

This section mainly concerns qualitative aspects of the generalized
Fermat-Torricelli problem \eqref{ft} related to deriving necessary
as well as necessary and sufficient conditions for its solutions in
convex and nonconvex cases. We also show that the obtained
qualitative results allow us to explicitly find generalized
Fermat-Torricelli points in some remarkable
settings.\vspace*{0.05in}

Let us first establish sufficient conditions for the {\em existence}
of optimal solutions to the generalized Fermat-Torricelli problem
under consideration.

\begin{Proposition}\label{exist} {\bf (existence of optimal
solutions to the generalized Fermat-Torricelli problem).} In
addition to the standing assumption of Section~{\rm 1}, suppose that
at least one of the sets $\O_1,\ldots,\O_n$ in \eqref{ft} is bounded
and that $\inf_{x\in X}T(x)<\infty$. Then the generalized
Fermat-Torricelli problem \eqref{ft} admits an optimal solution in
each of the following settings:

{\bf (i)} The space $X$ is finite-dimensional.

{\bf (ii)} The space $X$ is reflexive and all the sets $\O_i$,
$i=1,\ldots,n$, are convex.
\end{Proposition}
{\bf Proof.} To justify (i), suppose that the set $\O_1$ is bounded.
Denoting $\al:=\inf_{x\in X}T(x)$, we immediately observe that
\begin{equation*}
\big\{x\in X\big|\;T(x)<\al+1\big\}\subset\big\{x\in X\big|\;
T^F_{\O_1}(x)<\al+1\big\},
\end{equation*}
and hence the level set $\{x\in X|\; T(x)<\al+1\}$ is bounded. By
\cite[Proposition~3.5]{bmn10} the minimal time function
\eqref{minimal time} is lower semicontinuous under the assumptions
in (i). Thus we deduce the existence of solutions to \eqref{ft} from
the classical Weierstrass theorem.

To proceed with the proof of (ii), recall that every convex,
bounded, and closed subset of a reflexive space is sequentially
weakly compact. Furthermore, Proposition~3.5 of \cite{bmn10} yields
the lower semicontinuity of the minimal time function \eqref{minimal
time} with such a target set. This implies the weak lower
semicontinuity of $T(\cdot)$ from \eqref{ft} under the convexity
assumptions made and hence ensures the existence of optimal
solutions to \eqref{ft} in case (ii) by applying the Weierstrass
theorem in the weak topology of $X$ to this problem.
$\h$\vspace*{0.05in}

It is not hard to illustrate by examples that all the assumptions
made in Proposition~\ref{exist} are {\em essential} for the
existence of optimal solutions to \eqref{ft}. Consider for instance
a particular case of \eqref{distance function} with $X=\R^2$, $n=2$,
$\O_1:=\{(x,y)\in \R^2\;|\;y\ge e^x\}$, and $\O_2:=\R\times\{0\}$.
It is clear that this problem does not have an optimal
solution.\vspace*{0.05in}

Let us further proceed with deriving {\em optimality conditions} for
the generalized Fermat-Torricelli problem \eqref{ft}. Define the
sets
\begin{equation}\label{A}
A_i(x):=\bigcup_{\omega\in\Pi^F_{\O_i}(x)}\big[-\partial
\rho_F(\omega-x)\cap N(\omega;\O_i)\big],\quad x\in X,
\end{equation}
provided that $\Pi^F_{\O_i}(x)\ne\emp$. As in the proof of
Proposition~\ref{exist} (with no need of boundedness of the target
set $\O$ while under the standing assumptions on the dynamics $F$),
we can deduce from the generalized projection definition \eqref{pr}
that $\Pi^F_\O(x)\ne\emp$ for all $x\notin\O$ in each of the two
following cases:\vspace*{0.05in}

$\bullet$ $X$ is finite-dimensional and $\O$ is closed;

$\bullet$ $X$ is reflexive, $\O$ is closed and
convex.\vspace*{0.05in}\\
Furthermore, it easy to observe from the construction in \eqref{A}
that
\begin{eqnarray}\label{A1}
A_i(x)=N(x;\O)\cap C^*\;\mbox{ as }\;x\in\O_i,\quad i=1,\ldots,n,
\end{eqnarray}
for arbitrary closed sets $\O_i$, where the support level set $C^*$
is defined in \eqref{sl}. Useful relationships between the sets
$A_i(x)$ and the subdifferential $\partial T^F_{\O_i}(x)$ in the
out-of-set case $x\in\O_i$ follow from Theorem~\ref{mink3} and
Theorem~\ref{convex case 2} for convex and nonconvex targets $\O_i$.
These relationships are widely used in the sequel.\vspace*{0.05in}

We first establish necessary optimality conditions for the general
nonconvex problem \eqref{ft} in infinite dimensions. For simplicity
we assume that $0\in\i F$, which ensures the Lipschitz continuity of
the minimal time function $T^F_{\O_i}(\cdot)$ for all sets $\O_i$,
$i=1,\ldots,n$, in \eqref{ft} and the possibility to apply the sum
rule from Theorem~\ref{mink3}. Our approach allows us to treat, with
some elaboration, the non-Lipschitzian case when $\i F=\emp$ by
using more involved subdifferential formulas for the minimal time
function obtained in \cite{bmn10} and the basic subdifferential sum
rules for non-Lipschitzian functions given in
\cite[Chapter~3]{mor06a}.

\begin{Theorem}\label{necessary} {\bf (necessary optimality conditions
for the generalized Fermat-Torricelli problem).} Let $X$ be an
Asplund space, and let $0\in\i F$. If $\ox\in X$ is a local optimal
solution to the generalized Fermat-Torricelli problem \eqref{ft}
such that for each $i=1,\ldots,n$ the minimal time function
$T^F_{\O_i}(\cdot)$ is well posed at $\ox$ when $\ox\notin\O_i$,
then
\begin{equation}\label{A2}
0\in\sum_{i=1}^n A_i(\ox)
\end{equation}
with the sets $A_i(\ox)$, $i=1,\ldots,n$, defined in \eqref{A}.
\end{Theorem}
{\bf Proof.} If $\ox$ is a local solution to \eqref{ft}, then
$0\in\partial T(\ox)$ by the generalized Fermat stationary rule; see
\cite[Proposition~1.114]{mor06a}. It is well known that the minimal
time function \eqref{minimal time} is Lipschitz continuous on $X$
provided that $0\in\i F$; see, e.g., \cite[Lemma~3.2]{heng}.
Employing thus the nonconvex subdifferential sum rule from
Theorem~\ref{nonconvex sum rule}, we have
\begin{equation}\label{A4}
0\in\sum_{i=1}^n\partial T^F_{\O_i}(\ox).
\end{equation}
Comparing inclusion \eqref{limiting1} from Theorem~\ref{lim-min}
with formula \eqref{A1} gives us the upper estimate
\begin{equation}\label{A3}
\partial T^F_{\O_i}(\ox)\subset A_i(\ox),\quad i=1,\ldots,n,
\end{equation}
in the in-set case $\ox\in\O_i$. Furthermore, by Theorem~\ref{mink3}
the above inclusion \eqref{A3} holds also in the out-of-set case
$\ox\notin\O_i$ under the assumed well-posedness. Substituting
\eqref{A3} into \eqref{A4}, we arrive at \eqref{A2} and complete the
proof of the theorem. $\h$\vspace*{0.05in}

For the particular case \eqref{distance function} of problem
\eqref{ft} in Hilbert spaces a more explicit counterpart of
\eqref{A2} holds, which provides a convenient necessary optimality
condition for the Steiner-type extension \eqref{distance function}
of the Fermat-Torricelli problem.

\begin{Theorem}\label{nonconvex} {\bf (necessary optimality
conditions for the Steiner-type extension of the Fermat-Torricelli
problem in Hilbert spaces).} Let $X$ be a Hilbert space, and let
$\ox\in X$ be a local optimal solution to problem \eqref{distance
function} such for each $i=1,\ldots,n$ the distance function ${\rm
d}(\cdot;\O_i)$ is well posed at $\ox$ when $\ox\notin\O_i$. Then
condition \eqref{A2} is necessary for optimality of $\ox$ in
\eqref{distance function}, where the sets $A_i(\ox)$ are explicitly
expressed by:
\begin{eqnarray}\label{A5}
A_i(\ox)=\left\{\begin{array}{lr}
\dfrac{\ox-\Pi(\ox;\O_i)}{\rm{d}(\ox;\O_i)}&\mbox{ if
}\;\ox\notin\O_i\\\\
N(\ox;\O_i)\cap\B&\mbox{ if }\;\ox\in\O_i
\end{array}
\right.\;\mbox{ for all }\;i=1,\ldots,n.
\end{eqnarray}
\end{Theorem}
{\bf Proof.} Observe that $\sigma(\cdot)=\|\cdot\|$ in \eqref{sf}
and that $C^*=\B$ in \eqref{sl} by $X^*=X$ and $F=\B$. By
Theorem~\ref{necessary} it remains to prove that the expression for
$A_i(\ox)$ in \eqref{A5} for $\ox\notin\O_i$ reduces to \eqref{A} in
the setting under consideration.

Fix an arbitrary vector $\ox\notin\O_i$ and show that any vector
$u\in A_i(\ox)$ from \eqref{A} belongs to the set on the right-hand
side of \eqref{A5}. It is well known that in Hilbert spaces we have
\begin{equation*}
\la\ox-\omega,x-\omega\ra\le\frac{1}{2}\|x-\omega\|^2\;\mbox{ for
all }\;\omega\in\Pi(\ox;\O_i)\;\mbox{ and }\;x\in\O_i.
\end{equation*}
The latter implies, by definitions \eqref{1.2.1} and \eqref{2.3.1},
that
\begin{eqnarray*}
\ox-\omega\in\Hat N(\omega;\O_i)\subset N(\omega;\O_i).
\end{eqnarray*}
Since $\rho_F(x)=\|x\|$ for the Minkowski gauge \eqref{mk} in this
case, it gives
\begin{equation}\label{A6}
-\partial\rho(\omega-\ox)=\Big\{\frac{\ox-\omega}{\|\ox-\omega\|}\Big\}.
\end{equation}
Using now \eqref{A6} and the inclusion
\begin{equation*}
\frac{\ox-\omega}{\|\ox-\omega\|}\in N(\omega;\O_i)
\end{equation*}
held by representation \eqref{nc}, we have for $u\in A_i(\ox)$ in
\eqref{A} the relationships
\begin{equation*}
u\in-\partial\rho_F(\omega-\ox)\cap
N(\omega;\O_i)=\Big\{\frac{\ox-\omega}{\|\ox-\omega\|}\Big\}.
\end{equation*}
It follows that $\|\ox-\omega\|=\mbox{d}(\ox;\O_i)$ due to
$\omega\in\Pi(\ox;\O_i)$, and thus
\begin{equation*}
u=\frac{\ox-\omega}{\|\ox-\omega\|}\in\frac{\ox-\Pi(\ox;\O_i)}{\rm{d}(\ox;\O_i)},
\end{equation*}
which justifies that $u$ belongs to the set on the right-hand side
of \eqref{A5}.

To prove the converse inclusion, take any
$u\in\disp\frac{\ox-\Pi(\ox;\O_i)}{\rm{d}(\ox;\O_i)}$ and find
$\omega\in\Pi(\ox;\O_i)$ such that
$u=\disp\frac{\ox-\omega}{\|\ox-\omega\|}$. Then $u\in
N(\omega;\O_i)$ by \eqref{nc} and $u\in-\partial\rho_F(\omega-\ox)$
by \eqref{A6}. This gives
\begin{eqnarray*}
u\in-\partial\rho_F(\omega-\ox)\cap N(\omega;\O_i)
\end{eqnarray*}
and shows that the set $A_i(\ox)$ in \eqref{A5} belongs to the one
in \eqref{A} for each $i\in\{1,\ldots,n\}$, which thus completes the
proof of the theorem. $\h$\vspace*{0.05in}

Observe that the well-posedness assumption of
Theorem~\ref{nonconvex} is automatic when either $X$ is
finite-dimensional or the corresponding set $\O_i$ is convex.
Furthermore, we also have $\Pi(\ox;\O_i)\ne\emp$ in the same
settings.\vspace*{0.05in}

Next we employ Theorem~\ref{nonconvex} to specify optimal solutions
to \eqref{distance function} with $n=3$ therein. Note that the
condition $\la u,v\ra\le-1/2$ obtained in what follows means that
the angle between these two vectors is larger than or equal to
$120^\circ$, which is the crucial case in the classical
Fermat-Torricelli problem.

\begin{Corollary}\label{necessary conditions} {\bf (necessary conditions
for the generalized Fermat-Torricelli problem with three nonconvex
sets in Hilbert spaces).} Let $n=3$ in the framework of
Theorem~{\rm\ref{nonconvex}}, where $\O_1,\O_2,\O_3$ are pairwise
disjoint subsets of $X$. The following alternative holds for a local
optimal solution $\ox\in X$ with the sets $A_i(\ox)$ defined by
\eqref{A5}:

{\bf (i)} The point $\ox$ belongs to one of the sets $\O_i$, say
$\O_1$, and does not belong to the two others. Then there are
$a_2\in A_2(\ox)$ and $a_3\in A_3(\ox)$ such that
\begin{equation}\label{A6a}
\big\la a_2,a_3\ra\le-1/2\;\mbox{ and }\;-a_2-a_3\in N(\ox;\O_1).
\end{equation}

{\bf (ii)} The point $\ox$ does not belong to all the three sets
$\O_1$, $\O_2$, and $\O_3$. Then there are $a_i\in A_i(\ox)$ as
$i=1,2,3$ such that
\begin{equation}\label{A6b}
\la a_i,a_j\ra=-1/2\;\mbox{ for }\;i\ne j\;\mbox{ as
}\;i,j\in\big\{1,2,3\big\}.
\end{equation}
\end{Corollary}
{\bf Proof.} Since the sets $\O_i$ are pairwise disjoint, the
settings in (i) and (ii) fully describe the possible location of
$\ox$. Then Theorem~\ref{nonconvex} ensures that
\begin{equation}\label{A7}
0\in A_1(\ox)+A_2(\ox)+A_3(\ox)
\end{equation}
with the sets $A_i(\ox)$ defined by \eqref{A5}. Considering first
case (i), we get by \eqref{A5} and \eqref{A7} vectors $a_2\in
A_2(\ox)$ and $a_3\in A_3(\ox)$ satisfying the relationships
\begin{equation}\label{A7a}
\|a_2\|=\|a_3\|=1\;\mbox{ and }\;-a_2-a_3\in N(\ox;\O_1)\cap\B.
\end{equation}
Due to the obvious identities
\begin{equation*}
\|a_2+a_3\|^2=\|a_2\|^2+2\la a_2,a_3\ra+\|a_3\|^2=2+2\la a_2,a_3\ra,
\end{equation*}
the condition $\|a_2+a_3\|\in\B\Longleftrightarrow\|a_2+a_3\|^2\le
1$ is equivalent to $\la a_2,a_3\ra\le-1/2$. Thus the necessary
optimality condition \eqref{A7} can be equivalently rewritten in
form \eqref{A6a}, which completes the proof in case (i).

Considering next case (ii), we get by \eqref{A7} and \eqref{A5}
vectors $a_i\in A_i(\ox)$ for $i=1,2,3$ satisfying the relationships
\begin{equation*}
\|a_1\|=\|a_2\|=\|a_3\|=1\;\mbox{ and }\;a_1+a_2+a_3=0.
\end{equation*}
The latter implies that $a_1+a_2=-a_3$, and hence
\begin{equation*}
\la a_1,a_3\ra+\la a_2,a_3\ra=\la a_1+a_2,a_3\ra=-\la a_3,a_3\ra=-1.
\end{equation*}
Similarly we arrive at the conditions
\begin{equation*}
\la a_1,a_2\ra+\la a_1,a_3\ra=-1\;\mbox{ and }\;\la a_1,a_2\ra+\la
a_2,a_3\ra=-1.
\end{equation*}
Subtracting pairwisely the above equalities gives us
\begin{equation*}
\la a_1,a_2\ra=\la a_2,a_3\ra=\la a_1,a_3\ra=-1/2,
\end{equation*}
which can be written as \eqref{A6b} and thus completes the proof of
the corollary. $\h$\vspace*{0.05in}

The following example illustrates the application of
Corollary~\ref{necessary conditions} to a particular problem on the
plane with two convex and one nonconvex sets.

\begin{figure}[htb]
\begin{center}
\includegraphics[height=2.2in,width=1.8in]{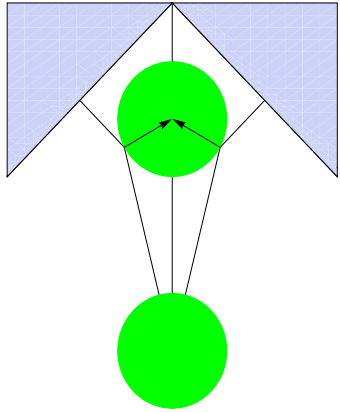}
\caption{A Nonconvex Generalized Fermat-Torricelli Problem.}
\end{center}
\end{figure}

\begin{Example}\label{ex-nonconvex} {\bf (nonconvex generalized
Fermat-Torricelli problem on the plane).} {\rm In the setting of
Corollary~\ref{necessary conditions}, let $\O_1$ be the ball
centered at $c_1:=(0,-2)$ with radius $r=1$, let $\O_2$ be the ball
centered at $c_2:=(0,-6)$ with the same radius $r=1$, and let $\O_3$
be a nonconvex set defined by
\begin{equation}
\O_3:=\big\{(x_1,x_2)\in\R^2\big|\;x_2\ge-|x_1|\big\}
\end{equation}
as depicted on Figure~1. By Proposition~\ref{exist} there is an
optimal solution to this problem. Applying then
Corollary~\ref{necessary conditions}, all the assumptions of which
are satisfied, we find two points lying on the boundary of $\O_1$
and denoted by $u$ and $v$ such that $[c_1,u]$ is the angle bisector
for the angle formed by the lines $uc_2$ and $up_u$, where $p_u$ is
the projection of $u$ to $\O_3$, while $[v,c_1]$ is the angle
bisector for the angle formed by the lines $vc_2$ and $vp_v$. These
two points satisfy {\bf(i)} of Corollary \ref{necessary conditions}
and they are actually the optimal solutions to the problem under
consideration. It is not hard to find $u$ and $v$ numerically; we
get $u=(-0.8706,-2.4920)$ and $v=(0.8706,-2.4920)$ up to five
significant digits, with the optimal value of the problem equal to
3.7609.}
\end{Example}

We continue now by considering the generalized Fermat-Torricelli
problem \eqref{ft} with {\em convex} target set $\O_i$ as
$i=1,\ldots,n$ in Banach spaces. In this case we derive {\em
necessary and sufficient} optimality conditions for
Fermat-Torricelli points.

It follows from Theorem~\ref{lim-min} and Theorem~\ref{convex case
2} that in the case of convex sets $\O_i$ with
$\Pi^F_{\O_i}(x;\O_i)\ne\emp$ as $x\notin\O_i$ and $0\in F$ we have
the equalities
\begin{equation}\label{A8}
A_i(x)=-\partial\rho_F(\omega-x)\cap N(\omega;\O_i)\;\mbox{ for any
}\;x\in X\;\mbox{ and }\;\omega\in\Pi^F_{\O_i}(x;\O_i)
\end{equation}
for the sets $A_i(x)$ defined in \eqref{A}, where the
subdifferential and normal cone are explicitly computed by formulas
\eqref{convex subdifferential} and \eqref{cnc} of convex analysis.
Here is a characterization of Fermat-Torricelli points for convex
problems.

\begin{Theorem}\label{necessary and sufficient convex} {\bf
(necessary and sufficient conditions for generalized
Fermat-Torricelli points of convex problems in Banach spaces).} Let
all the target sets $\O_i$ be convex, let $0\in\i F$ for problem
\eqref{ft} formulated in a Banach space $X$, and let $\ox\in X$ be
such that $\Pi^F_{\O_i}(\ox)\ne\emp$ whenever $\ox\notin\O_i$ as
$i=1,\ldots,n$. Then condition \eqref{A2} with the sets $A_i(\ox)$
defined in \eqref{A8} is necessary and sufficient for optimality of
$\ox$ in this problem.
\end{Theorem}
{\bf Proof.} As mentioned above, the convexity of all the sets
$\O_i$ implies the convexity of the cost function $T(x)$ in problem
\eqref{ft}. By the generalized Fermat rule \eqref{fermat} for convex
functions we get the inclusion $0\in\partial T(\ox)$ as a necessary
and sufficient conditions for optimality of $\ox\in X$ in
\eqref{ft}. Since all the functions $T^F_{\O_i}(\cdot)$ are locally
Lipschitzian under the interiority assumption on the dynamics $F$,
the convex subdifferential sum rule of Theorem~\ref{sum rule}
ensures that the latter inclusion is equivalent to
\begin{equation}\label{A8a}
0\in\sum_{i=1}^n\partial T^F_{\O_i}(\ox).
\end{equation}
Applying now relationship \eqref{A1} and equality \eqref{limiting2}
of Theorem~\ref{lim-min} in the in-set case as well as
Theorem~\ref{convex case 2} in the out-of-set case, we conclude that
$\partial T^F_{\O_i}(\ox)=A_i(\ox)$ as $i=1,\ldots,n$. Thus
inclusion \eqref{A8a} is equivalent to \eqref{A2}, and the latter is
necessary and sufficient for optimality of $\ox$ in the convex
Fermat-Torricelli problem \eqref{ft}. $\h$\vspace*{0.05in}

The following consequence of Theorem~\ref{necessary and sufficient
convex} provides an explicit characterization of Fermat-Torricelli
points in the convex Steiner-type extension \eqref{distance
function} of the classical problem in the Hilbert space setting. In
this case we use formula \eqref{A5} for constructing the sets
$A_i(\ox)$, which reduce to {\em singletons} if $\ox\notin\O_i$ and
are computed explicitly by \eqref{cnc} if $\ox\in\O_i$.

\begin{Corollary}\label{st-convex} {\bf (characterization of
optimal solutions to the convex Steiner-type extension of the
Fermat-Torricelli problem in Hilbert spaces).} Let $X$ be a Hilbert
space, and let all the sets $\O_i$ in \eqref{distance function} be
convex. Then condition \eqref{A2} with $A_i(\ox)$ computed in
\eqref{A5} is necessary and sufficient for optimality of $\ox\in X$
in problem \eqref{distance function}.
\end{Corollary}
{\bf Proof.} It follows from Theorem~\ref{necessary and sufficient
convex} due the fact that the sets $A_i(\ox)$ from \eqref{A} reduce
to those in \eqref{A5} as proved in Theorem~\ref{nonconvex} and due
to the projection nonemptiness $\Pi(\ox;\O_i)\ne\emp$ for any
$\ox\in X$ and $i=1,\ldots,n$ in the setting under consideration.
$\h$\vspace*{0.05in}

Note that, in contrast to problem \eqref{distance function}
addressed in Corollary~\ref{st-convex}, the characterization of
generalized Fermat-Torricelli points for problem \eqref{ft} obtained
in Theorem~\ref{necessary and sufficient convex} depends on the
dynamics $F$ and clearly determines different solutions for the same
targets sets $\O_i$ while different dynamics sets $F$. For example,
consider the case of the three singletons $\O_1:=\{(-1,0)\}$,
$\O_2:=\{(0,1)\}$, and $\O_3:=\{(1,0)\}$ on the plane $\R^2$. Then
Corollary~\ref{st-convex} gives us the unique optimal solution
$(0,1/\sqrt{3})$ to the corresponding problem \eqref{distance
function}, while for $F=[-1,1]\times[-1,1]$ and the same sets
$\O_i$, $i=1,2,3$, we have the unique optimal solution  $(0,1)$ to
the generalized Fermat-Torricelli problem \eqref{ft} with these sets
$F$ and $\O_i$.\vspace*{0.05in}

Let us present a simple application of Corollary~\ref{st-convex} to
a version of the generalized Fermat-Torricelli problem
\eqref{distance function} for finitely many disjoint closed
intervals of the real line.

\begin{Proposition}\label{line} {\bf (Fermat-Torricelli problem for
closed intervals of the real line).} Consider problem
\eqref{distance function} with the sets $\O_i$ given by $n$ disjoint
closed intervals $[a_i,b_i]\subset\R$ as $i=1,\ldots,n$, where
$a_1\le b_1<a_2\le b_2<\ldots<a_n\le b_n$. The following hold:

{\bf (i)} If $n=2k+1$, then any point of the interval
$[a_{k+1},b_{k+1}]$ $($the mid interval$)$ is an optimal solution to
the problem under consideration.

{\bf (ii)} If $n=2k$, then any point of the interval $[b_k,
a_{k+1}]$ is an optimal solution to the problem under consideration.
\end{Proposition}
{\bf Proof.} Let $\ph(x):=\mbox{d}(x;\O)$ with $\O=[a,b]$. It is
easy to compute (as, e.g., a particular case of Theorem~\ref{convex
case 2}) the subdifferential of $\ph$ by:
\begin{displaymath}
\partial\ph(\ox)=\left\{\begin{array}{ll}
\{0\} & \mbox{ if }\;a<\ox<b,\\
${\rm [-1,0]}$ & \mbox{ if }\;\ox=a,\\
${\rm [0,1]}$ & \mbox{ if }\;\ox=b,\\
${\rm \{-1\}}$ & \mbox{ if }\;\ox<a,\\
${\rm \{1\}}$ & \mbox{ if }\;\ox>b.
\end{array}
\right.
\end{displaymath}
Consider first case (i) when $n=2k+1$, we get for any $\ox\in
[a_{k+1},b_{k+1}]$ the relationships
\begin{align*}
&\sum_{i=1}^k\partial\mbox{d}(\ox;\O_i)=\{k\},\\
&\sum_{i=k+2}^n\partial\mbox{d}(\ox;\O_i)=\{-k\},\\
&\mbox{and } 0\in\partial\mbox{d}(\ox;\O_{k+1}).
\end{align*}
The latter implies that
$0\in\sum_{i=1}^n\partial\mbox{d}(\ox;\O_i)$, which ensures by
Corollary~\ref{st-convex} that $\ox$ is an optimal solution to the
problem under consideration. Taking further any
$\ox\notin[a_{k+1},b_{k+1}]$, we get by the above calculation that
$0\notin\sum_{i=1}^n\partial\mbox{d}(\ox;\O_i)$ and hence learn from
the characterization of Corollary~\ref{st-convex} that such a number
$\ox$ cannot be an optimal solution to the problem. The even case of
$n$ in (ii) is treated similarly. $\h$\vspace*{0.05in}

Another application of Corollary~\ref{st-convex} provides complete
characterizations of Fermat-Torricelli points for the convex problem
\eqref{distance function} with $n=3$ in Hilbert spaces. Note that in
this case, due the projection uniqueness, we have
\begin{equation}\label{a}
A_i(\ox)=a_i:=\Big\{\dfrac{\ox-\Pi(\ox;\O_i)}{\mbox{d}(\ox;\O_i)}
\Big\},\quad\ox\notin\O_i,
\end{equation}
for the sets $A_i(\ox)$ defined in \eqref{A5}.

\begin{Proposition}\label{necessary and sufficient conditions}
\label{3H} {\bf (characterizations of generalized Fermat-Torricelli
points for three convex sets in Hilbert spaces).} Let $X$ be a
Hilbert space, and let $\O_1,\O_2,\O_3$ be pairwise disjoint convex
subsets of $X$. Then $\ox\in X$ is an optimal solution to problem
\eqref{distance function} generated by these sets if and only if one
of the conditions {\rm (i)} and {\rm (ii)} of
Corollary~{\rm\ref{necessary conditions}} is satisfied, where the
vectors $a_i$, $i=1,2,3$, are defined in \eqref{a}, and where the
normal cone $N(\ox;\O_1)$ in \eqref{A6a} is computed by \eqref{cnc}.
\end{Proposition}
{\bf Proof.} The necessity part of the proposition follows from
Corollary~\ref{necessary conditions} by the observations that the
convex problems under consideration is well posed at $\ox$ and that
$\emp\ne\Pi(\ox;\O_i)$ is a singleton for any $i=1,2,3$. The
sufficiency part of the proposition can be derived from
Corollary~\ref{st-convex} by the arguments developed in the proof of
Corollary~\ref{necessary conditions}. $\h$\vspace*{0.05in}

Finally in this section, we illustrate the application of
Proposition~\ref{necessary and sufficient conditions} to some
particular problems of Fermat-Torricelli type formulated on the
plane.

\begin{figure}[htb]
\begin{center}
\includegraphics[height=1.8 in,width=2.3in]{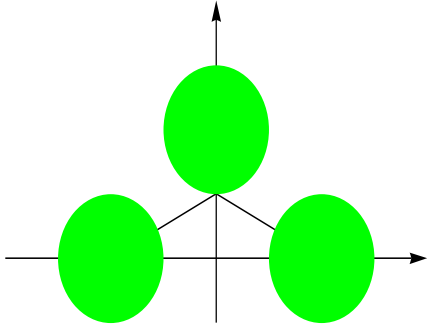}
\caption{A Convex Generalized Fermat-Torricelli Problem.}
\end{center}
\end{figure}

\begin{Example}\label{ex-convex} {\bf (convex generalized
Fermat-Torricelli problems on the plane).} {\rm Let the sets $\O_1$,
$\O_2$, and $\O_3$ in problem \eqref{distance function} are closed
balls in $\R^2$ of radius $r=1$ centered at the points $(0,2)$,
$(-2,0)$, and $(2,0)$, respectively; see Figure~2. We can easily see
that the point $(0,1)\in\O_1$ satisfies all the conditions in
Proposition~\ref{necessary and sufficient conditions}(i), and hence
it is an optimal solution (in fact a unique one) to this problem.

More generally, consider problem \eqref{distance function} in $\R^2$
generated by three arbitrary pairwise disjoint disks denoted by
$\O_i$, $i=1,2,3$. Let $c_1$, $c_2$, and $c_3$ be the centers of the
disks. Assume first that either the line segment $[c_2,c_3]\in\R^2$
intersects $\O_1$, or $[c_1,c_3]$ intersects $\O_2$, or $[c_1,c_2]$
intersects $\O_3$. It is not hard to check that any point of the
intersections (say of the sets $\O_1$ and $[c_2,c_3]$ for
definiteness) is an optimal solution to the problem under
consideration, since it satisfies the necessary and sufficient
optimality conditions of Proposition~\ref{necessary and sufficient
conditions}(i). Indeed, if $\ox$ is such a point, then $a_2$ and
$a_3$ from \eqref{a} are unit vectors with $\la a_2, a_3\ra=-1$ and
$-a_2-a_3=0\in N(\ox;\O_1)$.

If the above intersection assumptions are violated, we define three
points $q_1,q_2$, and $q_3$ as follows. Let $u$ and $v$ be the
intersections of $[c_1,c_2]$ and $[c_1,c_3]$ with the boundary of
the disk centered in $c_1$. Then we can see that there is a unique
point $q_1$ on the minor curve generated by $u$ and $v$ such that
the measures of angle $c_1q_1c_2$ and $c_1q_1c_3$ are equal. The
points $q_2$ and $q_3$ are defined similarly.
Proposition~\ref{necessary and sufficient conditions} yields that
whenever the angle $c_2q_1c_3$, or $c_1q_2c_3$, or $c_2q_3c_1$
equals or exceeds $120^\circ$ (say the angle $c_2q_1c_3$ does), then
the point $\ox:=q_1$ is an optimal solution to the problem under
consideration. Indeed, in this case $a_2$ and $a_3$ from \eqref{a}
are unit vectors with $\la a_2, a_3\ra\le-1/2$ and $-a_2-a_3\in
N(\ox;\O_1)$ because the vector $-a_2-a_3$ is orthogonal to $\O_1$.

If none of these angles equals or exceeds $120^\circ$, there is a
point $q$ not belonging to $\O_i$ as $i=1,2,3$ such that the angles
$c_1qc_2=c_2qc_3=c_3qc_1$ are of $120^\circ$, and $q$ is an optimal
solution to the problem. Observe that in this case the point $q$ is
also a unique optimal solution to the classical Fermat-Torricelli
problem determined by the points $c_1,c_2$, and $c_3$.}
\end{Example}

\section{Generalized Fermat-Torricelli Problem in Convex Settings:
Numerical Aspects}
\setcounter{equation}{0}

The concluding section of the paper is devoted to some numerical
aspects of solving the generalized Fermat-Torricelli problem
\eqref{ft} and its concretizations for the case of $n$ convex target
sets in finite-dimensional spaces. Based on the subgradient method
in convex optimization and the subdifferential calculus results
discussed in Sections~2 and 3, we develop a first-order algorithm of
solving a general convex problem \eqref{ft} and present some of its
specifications and implementations.

\begin{Theorem}\label{subgradient method2} {\bf (subgradient
algorithm for the generalized Fermat-Torricelli problem).} Let
$\O_i$, $i=1,\ldots,n$, be convex subsets of a finite-dimensional
Euclidean space $X$, let $0\in\i F$, and let $S\ne\emp$ be the set
of optimal solutions to problem \eqref{ft}. Picking a sequence
$\{\al_k\}$ as $k\in\N$ of positive numbers and a starting point
$x_1\in X$, consider the algorithm
\begin{equation}\label{al}
x_{k+1}=x_k-\al_k\sum_{i=1}^n v_{ik},\quad k=1,2,\ldots,
\end{equation}
with an arbitrary choice of vectors
\begin{equation}\label{a1}
v_{ik}\in-\partial\rho(\omega_{ik}-x_k)\cap
N(\omega_{ik};\O_i)\;\mbox{ for some
}\;\omega_{ik}\in\Pi^F_{\O_i}(x_k)\;\mbox{ if }\;x_k\notin\O_i
\end{equation}
and $v_{ik}=0$ otherwise. Assume that
\begin{equation}\label{a2}
\sum_{k=1}^\infty\alpha_k=\infty\;\mbox{ and }\;
\ell^2:=\sum_{k=1}^\infty\alpha_k^2<0.
\end{equation}
Then the iterative sequence $\{x_k\}$ in \eqref{a1} converges to an
optimal solution for problem \eqref{ft} and the value sequence
\begin{equation}\label{Vk}
V_k:=\min\big\{T(x_j)\big|\;j=1,\ldots,k\big\}
\end{equation}
converges to the optimal value $\Hat V$ in this problem.
Furthermore, we have the estimate
\begin{align*}
V_k-\Hat
V\le\dfrac{{\rm{d}}(x_1;S)^2+L^2\ell^2}{2\sum_{i=1}^k\alpha_k},
\end{align*}
where $0\le L<\infty$ is a Lipschitz constant of the function
$T(\cdot)$ from \eqref{ft} on $X$.
\end{Theorem}
{\bf Proof.} As mentioned above, the value function $T(\cdot)$ in
\eqref{ft} is convex and globally Lipschitzian on $X$. By
Theorem~\ref{lim-min} and Theorem~\ref{convex case 2} the convex
subdifferential of the minimal time functions \eqref{minimal time}
at $x_k$ is computed by
\begin{eqnarray}\label{al-sub}
\partial T^F_{\O_i}(x_k)=\left\{\begin{array}{ll}
N(x_k;\O_i)\cap\big\{v\in X\big|\;\sigma_F(-v)\le 1\big\}&\mbox{if
}\;x_k\in\O_i,\\\\
N(\omega_{ik};\O_i)\cap\big[-\partial\rho(\omega_{ik}-x_k)\big]&\mbox{
if }\;x_k\notin\O_i,
\end{array}\right.
\end{eqnarray}
where $\omega_{ik}\in\Pi^F_{\O_i}(x_k)$ is any generalized
projection vector, $i\in\{1,\ldots,n\}$, and $k\in\N$. Recalling now
the subgradient algorithm for minimizing the convex function
$T(\cdot)$, we have
\begin{equation}\label{al-sub1}
x_{k+1}=x_k-\al_k\sum_{i=1}^n v_k\;\mbox{ with }\;v_k\in\partial
T(x_k),\quad 1,2,\ldots.
\end{equation}
The convex subdifferential sum rule of Theorem~\ref{sum rule}
provides the representation
\begin{equation*}
v_k=\sum_{i=1}^nv_{ik}\;\mbox{ with }\;v_{ik}\in\partial
T^F_{\O_i}(x_k)
\end{equation*}
of the subgradient  $v_k$ in \eqref{al-sub1}. Substituting the
latter into \eqref{al-sub1} gives us algorithm \eqref{al} with
$v_{ik}$ satisfying \eqref{a1}. Employing now the well-known results
on the subgradient method for convex functions in the so-called
``square summable but not summable case" (see, e.g., \cite{bert}),
we arrive at the conclusions of the theorem under the conditions in
\eqref{a2}. $\h$\vspace*{0.05in}

Note that, using the above arguments, we can similarly apply to the
generalized Fermat-Torricelli problem the subgradient method for
convex optimization in the other cases considered in \cite{bert}
with the corresponding replacements of the convergence conditions
\eqref{a2}.\vspace*{0.05in}

Let us present a useful consequence of Theorem~\ref{subgradient
method2} in the setting of \eqref{ft} when the Minkowski gauge
\eqref{mk} is differentiable everywhere but the origin; this holds,
e.g., for the distance function \eqref{distance function1}. In the
case under consideration we denote
\begin{equation}\label{al0}
g_F(x):=\left\{\begin{array}{ll}\nabla\rho_F(x)&\mbox{if }\;x\ne
0,\\\\
0&\mbox{if }\;x=0.
\end{array}\right.
\end{equation}

\begin{Corollary}\label{subgradient method1} {\bf (subgradient
algorithm under smoothness assumptions).} In the setting of
Theorem~{\rm\ref{subgradient method2}}, assume in addition that the
Minkowski gauge $\rho_F(\cdot)$ is differentiable at every point
$X\setminus\{0\}$. Picking a sequence of positive numbers
$\{\al_k\}$ satisfying conditions \eqref{a2} and given a starting
point $x_1\in X$, form the algorithm
\begin{equation}\label{al1}
x_{k+1}=x_k +\alpha_k\sum_{i=1}^n g_F(\omega_{ki}-x_k),
\end{equation}
where $\omega_{ik}\in\Pi^F_{\O_i}(x_k)$ is an arbitrary projection
vector. Then all the conclusions of Theorem~{\rm\ref{subgradient
method2}} hold true for algorithm \eqref{al1}.
\end{Corollary}
{\bf Proof.} Fix $i\in\{1,\ldots,n\}$ and $k\in\N$. When
$x_k\notin\O_i$ we have $\omega_{ik}-x_k\ne 0$ for any
$\omega_{ik}\in\Pi^F_{\O_i}(x_k)$. Hence
$\partial\rho_F(w_{ik}-x_k)=\nabla\rho_F(x_k-\omega_{ik})=g_F(\omega_{ik}-x_k)$
by \eqref{al0}, and the intersection
$-\partial\rho_F(\omega_k-x_{ik})\cap N(\omega_{ik};\O_i)$ reduces
to the singleton $\{-g_F(\omega_{ik}-x_k)\}$, which we take for
$v_{ik}$ in Theorem~\ref{subgradient method2} when $x_k\notin\O_i$.
In the other hand, for $x_k\in\O_i$ we get
\begin{equation*}
v_{ik}:=0=-g_F(0)=-g_F(\omega_{ik}-x_k).
\end{equation*}
Thus algorithm \eqref{al1} in both cases agrees with \eqref{al}
under the assumptions made.$\h$\vspace*{0.05in}

We have a further specification of algorithm \eqref{al1} for the
convex problem \eqref{distance function}.

\begin{Corollary}\label{subg-dist} {\bf (subgradient algorithm for
convex Steiner-type extensions).} Consider problem \eqref{distance
function} with convex sets $\O_i$, $i=1,\ldots,n$, in a
finite-dimensional Euclidean space $X$. Given a sequence $\{\al_k\}$
of positive numbers satisfying \eqref{a2} and a starting point
$x_1\in X$, form algorithm \eqref{al1} with $g_F(\cdot)$ computed by
\begin{eqnarray}\label{al2}
g_F(\omega_{ki}-x_k)=\left\{\begin{array}{ll}
\dfrac{\Pi(x_k;\O_i)-x_k}{{\rm{d}}(x_k;\O_i)}&\mbox{ if
}\;x_k\notin\O_i,\\\\
0&\mbox{ if }\;x_k\in\O_i.
\end{array}\right.
\end{eqnarray}
Then all the conclusions of Theorem~{\rm\ref{subgradient method2}}
are satisfied for this algorithm.
\end{Corollary}
{\bf Proof.} Follows from Corollary~\ref{subgradient method1} with
$\rho_F(x)=\|x\|$ and $\nabla\rho_F(x)=\dfrac{x}{\|x\|}$ if $x\ne
0$. $\h$ \vspace*{0.05in}

Now we consider some examples of implementing the above subgradient
algorithms to the numerical solution of particular versions of the
generalized Fermat-Torricelli problem.

\begin{Example}\label{disk} {\bf (Fermat-Torricelli problem for
disks).} {\rm Consider the Steiner-type extension \eqref{distance
function} of the Fermat-Torricelli problem  for pairwise disjoint
circular disks in $\R^2$. Let $c_i=(a_i,b_i)$ and $r_i$,
$i=1,\ldots,n$, be the centers and the radii of the disks under
consideration. The subgradient algorithm of
Corollary~\ref{subg-dist} is written in this case as
\begin{equation}\label{al3}
x_{k+1}=x_k-\al_k\sum_{i=1}^n q_{ik},
\end{equation}
where the quantities $q_{ik}$ are given by
\begin{equation*}
q_{ik}=\left\{\begin{array}{ll} 0 &\mbox{if }\;\|x_k-c_i\|\le
r_i,\\\\
\disp\frac{x_k-c_i}{\|x_k-c_i\|} &\mbox{if }\;\|x_k-c_i\|>r_i.
\end{array}\right.
\end{equation*}
The corresponding quantities $V_k$ are evaluated by formula
\eqref{Vk} with
\begin{equation*}
T(x_j)=\sum_{i=1,\;x_j\notin\O_i}^n\big(\|x_j-c_i\|-r_i\big).
\end{equation*}

Writing a MATLAB program, we can compute by the above expressions
the values of $x_k$ and $V_k$ for any number of disks and
iterations. This allows us, in particular, to examine the
convergence of the algorithm in various settings. The following
table shows the results from the implementation of the above
algorithm for three circles with centers $(-2,0)$, $(0,2)$, and
$(2,0)$ and with the same radius $r=1$. The presented calculations
are performed for the sequence $\alpha_k=1/k$ satisfying \eqref{a2}
and the starting point $x_1=(5,7)$.
\begin{center}
\begin{tabular}{|c|c|c|}
\hline
\multicolumn{3}{|c|}{MATLAB RESULT}\\
\hline $k$ & $x_k$ & $V_k$\\
\hline 10 & (0.6224, 1.1995) & 2.7243\\
\hline 100 & (0.0552, 0.9984) & 2.4741\\
\hline 1,000 & (0.0047, 0.9995) & 2.4721\\
\hline 10,000 & (0.0004, 0.9999) & 2.4721\\
\hline 100,000 &(0.0000, 1.0000) & 2.4721\\
\hline 1,000,000 &(0.0000, 1.0000) & 2.4721\\
\hline
\end{tabular}
\end{center}
Observe that the numerical results obtained in this case are
consistent with the theoretical ones given in
Proposition~\ref{necessary and sufficient conditions}.

For four disks centered at $(0,0),(2,2),(1,0)$, and $(2,-2)$ and the
same radius $r=1/4$, the MATLAB program gives us the optimal point
$(0.8453,-0.0000)$ and the optimal value $4.7141$. For five disks
centered at $(-1,0),(-1,1),(0,2),(1,1)$, and $(1,0)$ with radius
$r=1/2$, we get the optimal solution $(0.0000, 0.8505)$ and the
optimal value equal to $3.2973$.}
\end{Example}

Next we apply the subgradient algorithm \eqref{al1} to the
Steiner-type extension \eqref{distance function} of the
Fermat-Torricelli problem for the case of squares $\O_i$, which is
significantly different from the case of disks in
Example~\ref{disk}.

\begin{Example}\label{square} {\bf (Fermat-Torricelli
problem for squares).} {\rm  Consider problem \eqref{distance
function} generated by pairwise disjoint squares $\O_i$,
$i=1,\ldots,n$, of {\em right position} in $\R^2$, i.e., such that
the sides of each square are parallel to the $x$-axis or the
$y$-axis; see Figure~3. The center of a square is the intersection
of its two diagonals, and its radius equals one half of the side.
Let $c_i=(a_i,b_i)$ and $r_i$, $i=1,\ldots,n$, be the centers and
the radii of the squares under considerations. Then the vertices of
the $i$th square are denoted by $v_{1i}=(a_i+r_i,b_i+r_i)$,
$v_{2i}=(a_i-r_i,b_i+r_i)$, $v_{3i}=(a_i-r_i,b_i-r_i)$, and
$v_{4i}=(a_i+r_i,b_i-r_i)$.

Given a starting point $x_1$ and a sequence $\{\al_k\}$ satisfying
the conditions in \eqref{a2}, the subgradient algorithm of
Corollary~\ref{subg-dist} can be written in form \eqref{al3}, where
$x_k=(x_{1k},x_{2k})$ and where the quantities $q_{ik}$ are computed
as follows:
\begin{eqnarray*}
q_{ik}=\left\{\begin{array}{ll} 0&\mbox{ if }\;|x_{1k}-a_i|\le
r_i\;\mbox{ and }\;|x_{2k}-b_i|\le
r_i,\\\\
\dfrac{x_k-v_{1i}}{\|x_k-v_{1i}\|}&\mbox{ if }\;x_{1k}-a_i>r_i\;
\mbox{ and }\;x_{2k}-b_i>r_i,\\\\
\dfrac{x_k-v_{2i}}{\|x_k-v_{2i}\|}&\mbox{ if
}\;x_{1k}-a_i<-r_i\;\mbox{ and }\;x_{2k}-b_i>r_i,\\\\
\dfrac{x_k-v_{3i}}{\|x_k-v_{3i}\|}&\mbox{ if }\;x_{1k}-a_i<-r_i\;
\mbox{ and }\;x_{2k}-b_i<-r_i,\\\\
\dfrac{x_k-v_{3i}}{\|x_k-v_{4i}\|}&\mbox{ if }\;x_{1k}-a_i>r_i\;
\mbox{ and }\;x_{2k}-b_i<-r_i,\\\\
(0,1) & \mbox{ if }\;|x_{1k}-a_i|\le r_i\;\mbox{ and
}\;x_{2k}-b_i>r_i,\\\\
(0,-1) & \mbox{ if }\;|x_{1k}-a_i|\le r_i\;\mbox{ and }
\;x_{2k}-b_i<-r_i,\\\\
(1,0) & \mbox{ if } \;x_{1k}-a_i>r_i\;\mbox{ and }\;|x_{2k}-b_i|\le
r_i,\\\\
(-1,0) & \mbox{ if }\;x_{1k}-a_i<-r_i\;\mbox{ and }\;|x_{2k}-b_i|\le
r_i
\end{array}\right.
\end{eqnarray*}
for all $i\in\{1,\ldots,n\}$ and $k\in\N$, with the corresponding
value sequence $V_k$ defined by \eqref{Vk}.\vspace*{0.05in}

Considering the implementation of the above algorithm for three
squares with centers $(-2,0)$, $(0,2)$, and $(2,0)$ and with the
same radius $r=1/2$, we arrive at the optimal solution $(0,1.3660)$
and the optimal value $3.5981$. At the same time applying the
theoretical results of Corollary~\ref{subg-dist} to this case gives
us the exact optimal solution $(0,\dfrac{\sqrt{3}+1}{2})$ with the
optimal value $\dfrac{2+3\sqrt{3}}{2}$, which are consistent with
the numerical calculation.

\begin{figure}[htb]
\begin{center}
\includegraphics[height=2.2in,width=2.2in]{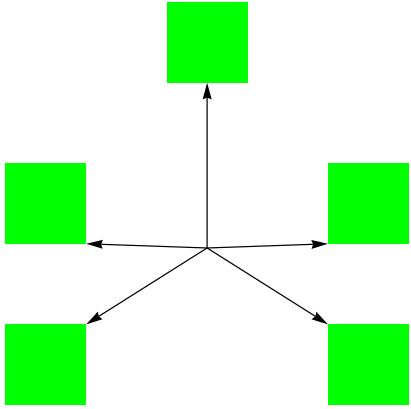}
\caption{A Generalized Fermat-Torricelli Problem for Five Squares.}
\end{center}
\end{figure}

For five squares with centers at $(-1,0),(-1,1),(0,2),(1,1)$, and
$(1,0)$ and the same radius $r=1/4$ we get by the subgradient
algorithm the optimal solution $(0.0000,0.7242)$ with the optimal
value equal to $4.3014$. However, it does not seem to be an easy
exercise to solve the above problem theoretically by using
Corollary~\ref{subg-dist}.}
\end{Example}

Let us finally illustrate applications of the subgradient algorithm
of Theorem~\ref{subgradient method2} to solving the generalized
Fermat-Torricelli problem \eqref{ft} formulated via the minimal time
function \eqref{minimal time} with non-ball dynamics. For
definiteness we consider the dynamics $F$ given by the square
$[-1,1]\times [-1,1]$ on the plane. In this case the corresponding
Minkowski gauge \eqref{mk} is given by the formula
\begin{equation}\label{mk1}
\rho_F(x_1,x_2)=\max\big\{|x_1|,|x_2|\big\}.
\end{equation}

Note that the function $\rho_F(\cdot)$ fails to be differentiable at
every nonzero point of $\R^2$, so we have to relay on the
subgradient algorithm of Theorem~\ref{subgradient method2} but not
of its corollaries above. Observe also that to implement algorithm
\eqref{al} we need to know just one element $v_{ik}$ from the set on
the right-hand side of \eqref{a1} for each $i\in\{1,\ldots,n\}$ and
$k\in\N$. By Theorem~\ref{convex case 2} the latter set agrees with
the subdifferential of the minimal time function $T^F_{\O_i}(x_k)$.
\vspace*{0.05in}

In the following proposition we compute a subgradient of the minimal
time function \eqref{minimal time} generated by the Minkowski gauge
\eqref{mk1} and a square target in $\R^2$, which is used then to
construct a subgradient algorithm to solve the corresponding
Fermat-Torricelli problem.

\begin{Proposition}\label{mk2} {\bf (subgradients of minimal
time functions with square dynamics and targets).} Let
$F=[-1,1]\times[-1,1]$, and let $\O$ be a square of right position
in $\R^2$ centered at $c=(a,b)$ with radius $r>0$. Then a
subgradient $v(\ox_1,\ox_2)\in\partial T^F_\O(\ox_1,\ox_2)$(not
necessarily uniquely defined) of the minimal time function
$T^F_\O(x_1,x_2)$ at $(\ox_1,\ox_2)$ is computed by
\begin{eqnarray}\label{v}
v(\ox_1,\ox_2)=\left\{\begin{array}{ll}
(1,0)&\mbox{if }\;|\ox_2-b|\le\ox_1-a,\;\ox_1>a+r,\\\\
(-1,0)&\mbox{if }\;|\ox_2-b|\le a-\ox_1,\;\ox_1<a-r.\\\\
(0,1)&\mbox{if }\;|\ox_1-a|\le\ox_2-b,\;\ox_2>b+r,\\\\
(0,-1)&\mbox{if }\;|\ox_1-a|\le b-\ox_2,\;\ox_2<b-r,\\\\
0&\mbox{if }\;(\ox_1,\ox_2)\in\O.
\end{array}\right.
\end{eqnarray}
\end{Proposition}
{\bf Proof.} It is ease to see that the minimal time function
\eqref{minimal time} admits the representation
\begin{equation*}
T^F_\O(x)=\inf_{\omega\in\O}\rho_F(x;\O).
\end{equation*}
This implies by \eqref{mk1} and the above structure of $\O$ that
\begin{eqnarray*}
T^F_{\O}(x_1,x_2)=\left\{\begin{array}{ll}
x_1-(a+r)&\mbox{if }\;|x_2-b|\le x_1-a,\;x_1>a+r,\\\\
(a-r)-x_1&\mbox{if }\;|x_2-b|\le a-x_1,\;x_1<a-r,\\\\
x_2-(b+r)&\mbox{if }\;|x_1-a|\le x_2-b,\;x_2>b+r,\\\\
(b-r)-x_2&\mbox{if }\;|x_1-a|\le b-x_2,\;x_2< b-r,\\\\
0&\mbox{if }\;(x_1,x_2)\in\O
\end{array}
\right.
\end{eqnarray*}
for all $(x_1,x_2)\in\R^2$. Applying to \eqref{mk1} the well-known
subdifferential formula for maximum functions in convex analysis
allows us to compute

{\small \begin{eqnarray*}
\partial\rho_F(\ox_1,\ox_2)=\left\{\begin{array}{ll}
\big\{(v_1,v_2)\in\R^2\;\big|\;|v_1|+|v_2|\le 1\big\}&\mbox{if }\;
(\ox_1,\ox_2)=(0,0),\\\\
\{(0,1)\}&\mbox{if }\;|\ox_1|<\ox_2,\\\\
\{(0,-1)\}&\mbox{if }\;\ox_2<-|\ox_1|,\\\\
\{(1,0)\}&\mbox{if }\;\ox_1>|\ox_2|,\\\\
\{(-1,0)\}&\mbox{if }\;\ox_1<-|\ox_2|,\\\\
\big\{(v_1,v_2)\in\R^2\;\big|\;|v_1|+|v_2|=1,\;v_1\ge 0,\;
v_2\ge 0\big\}&\mbox{if }\;\ox_1=\ox_2>0,\\\\
\big\{(v_1,v_2)\in\R^2\;\big|\;|v_1|+|v_2|=1,\;v_1\ge 0,\;v_2\le
0\big\}&\mbox{if }\;\ox_1=-\ox_2>0,\\\\
\big\{(v_1,v_2)\in\R^2\;\big|\;|v_1|+|v_2|=1,\;v_1\le 0,\;v_2\le
0\big\}&\mbox{if }\;\ox_1=\ox_2<0,\\\\
\big\{(v_1,v_2)\in\R^2\;\big|\;|v_1|+|v_2|=1,\;v_1\le 0,\;v_2\ge
0\big\}&\mbox{if }\;\ox_1=-\ox_2<0.
\end{array}\right.
\end{eqnarray*}}

In this way we can check by Theorem~\ref{convex case 2} that the
vector $v(\ox_1,\ox_2)$ is a subgradient of $T(\cdot,\cdot)$ at
$(\ox,\ox_2)$, which completes the proof of the proposition.
$\h$\vspace*{0.05in}

Now we are able to implement the subgradient algorithm of
Theorem~\ref{subgradient method2} to the problem under
consideration.

\begin{Example}\label{square1} {\bf (implementation of the
subgradient algorithm).} {\rm Consider the generalized
Fermat-Torricelli problem \eqref{ft} with the dynamics
$F=[-1,1]\times[-1,1]$ and the square targets $\O_i$ of right
position centered at $(a_i,b_i)$ with radii $r_i$ as $i=1,\ldots,n$.
Given a sequence of positive numbers $\{\al_k\}$ satisfying
\eqref{a2} and a starting point $x_1$, construct the subgradient
algorithm \eqref{al} for the iterations $x_k=(x_{1k},x_{2k})$ in
Theorem~\ref{subgradient method2}, where the vectors $v_{ik}$ are
computed by Proposition~\ref{mk2} as
\begin{eqnarray*}
v_{ik}=\left\{\begin{array}{ll} (1,0)&\mbox{if }\;|x_{2k}-b_i|\le
x_{1k}-a_i\;\mbox{ and }\;x_{1k}>a_i+r_i,\\\\
(-1,0)&\mbox{if }\;|x_{2k}-b_i|\le a_i-x_{1k}\;\mbox{ and }\;
x_{1k}<a_i-r_i,\\\\
(0,1)&\mbox{if }\;|x_{1k}-a_i|\le x_{2k}-b_i\;\mbox{ and }\;
x_{2k}>b_i+r_i,\\\\
(0,-1)&\mbox{if }\;|x_{1k}-a_i|\le b_i-x_{2k}\;\mbox{ and }\;
x_{2k}<b_i-r_i,\\\\
(0,0)&\mbox{otherwise.}
\end{array}\right.
\end{eqnarray*}

Implementing this algorithm for the case of three squares centered
at $(-2,0)$, $(0,2)$, and $(2,0)$ with radius $r=1/2$,
$\alpha_k=1/k$ and the initial point $(1,1)$, we arrive at an
optimal solution $(0.0000,1.5000)$ and the optimal value equal to
$3.0000$. For five squares centered at $(-1,0), (-1,1), (0,2),
(1,1)$, and $(1,0)$ with radius $r=1/4$, we have the optimal
solution $(0.0000,1.0000)$ and the optimal value equal to $3.7500$.}
\end{Example}

{\bf Acknowledgements.} The authors are grateful to Jon Borwein,
Michael Overton, and Boris Polyak for their valuable comments and
references.


\small

\begin{thebibliography}{99}
\markboth{\centerline{\sc }}{\centerline{\sc B. Mordukhovich and N.
M. Nam}}

\bibitem{acco} Andersen, K.D., Christiansen, A., Conn, A.R.,
Overton, M.L.: An efficient primal-dual interior-point method for
minimizing a sum of Euclidean norms. SIAM J. Scient. Comp. {\bf 22},
243--262 (2000).

\bibitem{bert} Bertsekas, D., Nedic, A., Ozdaglar, A.: Convex Analysis and Optimization.
Athena Scientific, Boston (2003)

\bibitem{b} Bhattacharya, B.B.: On the Fermat-Weber point of a polynomial chain. ArHiv:
1004.2958v1 (2010)

\bibitem{bms} Boltynski, V., Martini, H., Soltan, V.: Geometric Methods and Optimization Problems.
Kluwer Academic, Dordrecht (1999)

\bibitem{bv} Borwein, J.M., Vanderwerff, J.D.: Convex Functions: Characterizations,
Constructions and Counterexamples. Cambridge University Press (2009)

\bibitem{bz} Borwein, J.M., Zhu, Q.J.: Techniques of Variational Analysis. Springer,
CMS Books in Mathematics {\bf 20}, Springer, New York (2005)

\bibitem{cgm} Colombo, G., Goncharov, V.V., Mordukhovich, B.S.: Well-posedness of minimal
time problem with constant dynamics in Banach spaces. Set-Valued
Var. Anal., to appear (2010)

\bibitem{cowo} Colombo, G., Wolenski, P.R.: The subgradient formula for the minimal time
function in the case of constant dynamics in Hilbert space. J.
Global Optim. {\bf 28}, 269–-282 (2004)

\bibitem{heng} He, Y., Ng, K.F.: Subdifferentials of a minimum time function in Banach spaces.
J. Math. Anal. Appl. {\bf  321}, 896--910 (2006)

\bibitem{k} Kuhn, H.W.: Steiner's problem revisited. Studies
Math. {\bf 10}, 52--70 (1974)

\bibitem{msw} Martini, H., Swanepoel, K.J., Weiss, G.: The Fermat-Torricelli problem in
normed planes and spaces. J. Optim. Theory Appl. {\bf 115}, 283--314
(2002)

\bibitem{mor76} Mordukhovich, B.S.: Maximum principle in problems of time optimal control
with nonsmooth constraints. J. Appl. Math. Mech. {\bf 40}, 960--969
(1976)

\bibitem{mor06a} Mordukhovich, B.S.: Variational Analysis and Generalized Differentiation,
I: Basic Theory. Grundlehren Series (Fundamental Principles of
Mathematical Sciences) {\bf 330}, Springer, Berlin (2006)

\bibitem{mor06b} Mordukhovich, B.S.: Variational Analysis and Generalized Differentiation, II:
Applications. Grundlehren Series (Fundamental Principles of
Mathematical Sciences) {\bf 331}, Springer, Berlin (2006)

\bibitem{mn10} Mordukhovich, B.S., Nam, N.M.: Limiting subgradients of minimal time
functions in Banach spaces. J. Global Optim. {\bf 46}, 615--633
(2010)

\bibitem{bmn10} Mordukhovich, B.S, Nam, N.M.: Subgradients of minimal time functions
under minimal assumptions. J. Convex Anal., to appear (2010)

\bibitem{p} Phelps, R.R.: Convex Functions, Monotone Operators and Differentiability,
2nd edition. Lecture Notes Math. {\bf 1364}, Springer, Berlin (1993)

\bibitem{rw} Rockafellar, R.R., Wets, R.J-B.: Variational Analysis. Grundlehren Series
(Fundamental Principles of Mathematical Sciences) {\bf 317},
Springer, Berlin (1998)

\bibitem{s} Schirotzek, W.: Nonsmooth Analysis. Universitext, Springer, Berlin (2007)

\bibitem{tan} Tan, T.V.: An extension of the Fermat-Torricelli problem. J. Optim.
Theory Appl., published online (2010)

\bibitem{wp} Weiszfeld, E.: On the point for which
the sum of the distances to $n$ given points is minimum. Ann. Oper.
Res. {\bf 167}, 7--41 (2009). Translated from the French original
[Tohoku Math. J. {\bf 43}, 335--386 (1937)] and annotated by Frank
Plastria.

\end{thebibliography}
\end{document}